\documentstyle[pb-diagram]{article}

\newcommand{\es}{\emptyset}

\newcommand{\ba}{\begin{array}}
\newcommand{\ea}{\end{array}}

\newtheorem{lemma}{LEMMA}[section]
\newtheorem{corollary}{COROLLARY}[section]
\newtheorem{definition}{DEFINITION}[section]
\newtheorem{proposition}{PROPOSITION}[section]

\newcommand{\be}{\begin{enumerate}}
\newcommand{\ee}{\end{enumerate}}
\newcommand{\bi}{\begin{itemize}}
\newcommand{\ei}{\end{itemize}}
\newcommand{\bd}{\begin{description}}
\newcommand{\ed}{\end{description}}

\newcommand{\beq}{\begin{eqnarray*}}
\newcommand{\eeq}{\end{eqnarray*}}
\newcommand{\seq}{\Rightarrow}

\author{ {F.Parlamento, F.Previale }
\\Department of Mathematics and Computer Science
\\University of Udine,  via  Delle Scienze 206, 33100 Udine, Italy.
\\Department of Mathematics
\\University of Turin, via Carlo Alberto 10, 10123 Torino, Italy
\\e-mail: {\em franco.parlamento$@$uniud.it},  {\em flavio.previale$@$unito.it}
}
 \title{Contraction Elimination in   Sequent Based Ground Equational Calculus }
\date{}

\begin{document}
\maketitle




\begin{abstract}

In   \cite{PP13} we have shown that the cut rule is eliminable in two  ground equational sequent calculi to be  denoted by $EQ_M$ and $EQ'$.  In this note we prove that the contraction rule is not eliminable in $EQ_M$ but it is eliminable in $EQ'$.

\end{abstract}

\subsection{$EQ_M$ and $EQ$}

\begin{definition}
$EQ_M$ is the calculus acting on  sequents with one formula in the succedent  having {\em atomic logical axioms} of the form $A\seq A$,   the reflexivity axioms  
$\seq t=t$ ($\seq =$);  the {\em left structural rules of weakening, exchange and contraction}:

\[
\ba{ccccc}
\Gamma\seq H&~~~&\Gamma_1,F,G,\Gamma_2\seq H&~~~&\Gamma,F, F\seq H\\
\cline{1-1}\cline{3-3}\cline{5-5}
\Gamma, F\seq D&&\Gamma_1, G, F, \Gamma_2\seq H&~~~&\Gamma,F\seq H
\ea
\]

 the {\em atomic cut rule}:
\[
\ba{c}
\Gamma \seq A~~~~\Lambda, A \seq H\\
\cline{1-1}
\Gamma,\Lambda  \seq H
\ea
\]

 and the {\em atomic   equality left  introduction 
rules} $=_1\seq$ and $=_2\seq$, namely:

\[
\ba{ccc}
\Gamma\seq D\{v/r\}&~~~~~~~&\Gamma\seq D\{v/r\}\\
\cline{1-1}
\cline{3-3}
\Gamma, r=s \seq D\{v/s\}&~~~~~~~&\Gamma, s=r\seq D\{v/s\}\
\ea
\]
where by {\em atomic} we mean that $A$ and $D$ are required to be atomic formulae.

\end{definition} 

{\bf Notation} In the following $A$ and $D$ will denote atomic formulae.

\

\begin{definition}
$EQ$ is the calculus acting on  sequents with one formula in the succedent  having  logical axioms of the form $A\seq A$, the reflexivity axioms  
$\seq t=t$ ($\seq =$);  the left structural rules of weakening, exchange and contraction, the atomic cut rule and the following atomic 
 {\em congruence rule} $CNG$:

\[
\ba{c}
\Gamma\seq D\{v/r\}~~~~~~~\Lambda \seq r=s\\
\cline{1-1}
\Gamma, \Lambda\seq D\{v/s\}
\ea
\]

\end{definition}

\begin{definition}
$cf.EQ$ and $cf.EQ_M$ denote the systems $EQ$ and $EQ_M$ deprived of the cut rule.

\end{definition}

\begin{proposition}\label{equiv}

EQ and $EQ_M$ are equivalent over the structural rules of weakening, exchange and cut, more precisely  the rules $=_1\seq$ and $=_2\seq$ are derivable in $EQ$, without using the contraction and the cut rule and, conversely, $CNG$ is derivable in $EQ_M$ without using the contraction rule.

\end{proposition}

{\bf Proof}
\[
\ba{ccc}
\Gamma\seq D\{v/r\}&~~~&r=s\seq r=s\\
\cline{1-3}

\multicolumn{3}{c}{\Gamma, r=s\seq D\{v/s\}}

\ea
\]

\[
\ba{ccc}
&&\seq s=s~~~~s=r\seq s=r\\
\cline{3-3}
\Gamma\seq D\{v/r\}&&s=r\seq r=s\\
\cline{1-3}
\multicolumn{3}{c}{\Gamma, s=r\seq D\{v/s\}}\
\ea
\]

\[
\ba{ccc}
&~~~&\Gamma \seq D\{v/r\}\\
\cline{3-3}
\Lambda \seq r=s &&\Gamma, r=s \seq D\{v/s\}\\
\cline{1-3}
\multicolumn{3}{c}{\Gamma, \Lambda \seq D\{v/s\}}
\ea
\]
$\Box$

\subsection{ Eliminating the Contraction Rule}

The contraction rule is not eliminable from  $EQ_M$. For example the sequent 

$a=f(a)\seq a=f(f(a))$, where $a$ is an individual parameter,  in $cf.EQ_M$  has the following derivation:

\[
\ba{c}
a=f(a)\seq a=f(a)\\
\cline{1-1}
a=f(a),a=f(a)\seq a=f(f(a))\\
\cline{1-1}
a=f(a)\seq a=f(f(a))
\ea
\]
but we can show that there is  no  derivation in $EQ_M$ of  $a=f(a)\seq a=f(f(a))$, that does not use the contraction rule.

\begin{definition} 
Let $EQ_M^-$ be obtained from $EQ_M$ by suppressing the contraction rule.
\end{definition}

{\bf Notation}  $\Gamma_=$ denotes the sequence  that is  obtained from $\Gamma$ by suppressing all the formulae that are not equalities.

\begin{proposition}
$\Gamma\seq r=s$ is  derivable in $EQ_M$ ($EQ_M^-$), if and only if  $\Gamma_=\seq r=s$ is derivable in 
$EQ_M$ ($EQ_M^-$) with a derivation that contains only equalities.

\end{proposition}

{\bf Proof} The "if" direction is immediate by  the weakening rule.

  The "only if " direction is established by induction on the height of a given derivation ${\cal D}$ of  $\Gamma\seq r=s$. If $h({\cal D})=0$ then either $\Gamma =\es$ and $r\equiv s$ or $\Gamma$ reduces to $r=s$. In both case the conclusion is obvious.
If $h({\cal D})>0$ and the last rule of ${\cal D}$ is not a cut, then, if a  principal formula is not an equality,  the induction hypothesis yields directly the desired derivation of $\Gamma_=\seq r=s$. Otherwise it suffices to apply the same rule to the derivation provided by the induction hypothesis.
Finally suppose ${\cal D}$ ha the form:

\[
\ba{ccc}
{\cal D}_0&~~~&{\cal D}_1\\
\Gamma\seq A&&\Lambda, A\seq r=s\\
\cline{1-3}\multicolumn{3}{c}{\Gamma,\Lambda\seq r=s}
\ea
\]
If $A$ is not an equality, by induction hypothesis applied to 
${\cal D}_1$ we have a derivation of $\Lambda_{=}\seq r=s$,
 from which we can obtain the desired derivation of $\Gamma_{=},\Lambda_{=}\seq r=s$ 
by weakenings and exchanges.  On the other hand if $A$ is $p=q$,
 by induction hypothesis applied to both ${\cal D}_0$ and ${\cal D}_1$, 
we have derivations of $\Gamma_{=}\seq p=q$ and $\Lambda_{=}, p=q\seq r=s$,
 from which the desired derivation of $\Gamma_=,\Lambda_=\seq r=s$ 
can be obtained by appying the cut rule. $\Box$

\

{\bf Notation}  $r\equiv s$ denotes that $r$ and $s$ are syntactically identical.

\begin{proposition}
 If $p_1=p_1, \ldots, p_n=p_n \seq r=s$ is derivable in 
$EQ_M$ ($EQ_M^-$), then $r\equiv s$, in particular if
 $\seq r=s$ is derivable in 
 $EQ_M$ ($EQ_M^-$), then $r\equiv s$.
\end{proposition}

{\bf Proof} By induction on the height of a given derivation ${\cal D}$ of 

$p_1=p_1, \ldots, p_n=p_n \seq r=s$ that, by the previous Proposition,  we may assume to involve  only equalities.
 The base case is immediate. The inductive step is also straightforward. Let us  deal with the case in which ${\cal D}$ ends with a cut, namely it has the form:

\[
\ba{ccc}
{\cal D}_0&~~~&{\cal D}_1\\
p_1=p_1,\ldots, p_i=p_i\seq p=q&&p_{i+1}=p_{i+1},\ldots, p_n=p_n, p=q \seq r=s\\
\cline{1-3}
\multicolumn{3}{c}{p_1=p_1,\ldots, p_n=p_n \seq r=s}
\ea
\]

By induction hypothesis applied to ${\cal D}_0$ we have that $p\equiv q$, then, by induction hypothesis applied to ${\cal D}_1$ (with $n+1-i$ in place of $n$) we conclude that $r\equiv s$. $\Box$

\begin{proposition}

If $E$ is an equality, then the following hold:

\bi
\item[a)] If $*)~~~p_1=p_1, \ldots, p_j=p_j, E, p_{j+1}=p_{j+1}, \ldots, p_n=p_n \seq a=f(f(a))$ is derivable in $EQ_M^-$, then $E$ coincides with $a=f(f(a))$ or with $f(f(a))=a$

\item[b)] If $**)~~~p_1=p_1, \ldots, p_j=p_j, E, p_{j+1}=p_{j+1}, \ldots, p_n=p_n \seq f(f(a))=a$ is derivable in $EQ_M^-$, then $E$ coincides with $a=f(f(a))$ or with $f(f(a))=a$
\ei
\end{proposition}

{\bf Proof} $a)$ and $b$ are proved symultaneously by induction on the height of derivations.

$a)$ Let ${\cal D}$ be a derivation  in $EQ_M^-$ of  $*)$. If $h({\cal D})=0$, then $n=0$ and $E$ coincides with $a=f(f(a))$. As for the induction step, let us first observe  that ${\cal D}$ cannot end with a weakening that introduces $E$, since, by the previous Proposition, $p_1=p_1,\ldots, p_n=p_n\seq a=f(f(a))$ is not derivable.
If ${\cal D}$ ends with an exchange the conclusion is immediate by  the induction hypothesis.
If ${\cal D}$ ends with a $=_1\seq $-inference, namely it has the form
\[
\ba{c}
{\cal D}_0\\
p_1=p_1,\ldots, p_n=p_n \seq r=s\\
\cline{1-1}
p_1=p_1,\ldots, p_n=p_n, E \seq a=f(f(a))
\ea
\]
then, by the previous Proposition, $r\equiv s$. The only possibilities  of obtaining $a=f(f(a))$ by a substitution applied  to $r=r$ is that $r\equiv f(f(a))$ or $r\equiv a$, in which case  $E$ is either $f(f(a))=a$ or $a=f(f(a))$. Similarly  if ${\cal D}$ ends with a $=_2\seq$-inference we have that $E$ is $a=f(f(a))$ or $f(f(a))=a$.

If ${\cal D}$ ends with a cut, we have two cases.

Case 1. (Assuming for notational  simplicity that $j=n$) ${\cal D}$ has the form:

\[
\ba{ccc}
{\cal D}_0&~~~~&{\cal D}_1\\
p_1=p_1, \ldots, p_i=p_i\seq A&&p_{i+1}=p_{i+1}, \ldots, p_n=p_n, E, A \seq a=f(f(a))\\
\cline{1-3}
\multicolumn{3}{c}{p_1=p_1, \ldots, p_n=p_n, E\seq a=f(f(a))}
\ea
\]

By the previous Proposition,  $A$ must be an identity, hence, by induction hypothesis (with $n$ replaced by $n+1$) applied to ${\cal D}_1$, $E$ must coincide either with $a=f(f(a))$ or with $f(f(a))=a$.

Case 2. ${\cal D}$ has the form:

\[
\ba{ccc}
{\cal D}_0&~~~~&{\cal D}_1\\
p_1=p_1,\ldots, p_j=p_j, E \seq A&&p_{j+1}=p_{j+1}, \ldots, p_n=p_n, A\seq a=f(f(a))\\
\cline{1-3}
\multicolumn{3}{c}{p_1=p_1, \ldots, p_n=p_n, E \seq a=f(f(a))}
\ea
\]

By induction hypothesis applied to ${\cal D}_1$, $A$ coincides  with $a=f(f(a))$ or with $f(f(a))=a$.
We can then apply the induction hypothesis, either case $a)$ or case $b)$, to ${\cal D}_0$, to conclude that
$E$ coincides with $a=f(f(a)$ or with $f(f(a))=a$. 

The proof of $b)$ is entirely similar. $\Box$

\ 

Thus $a=f(a) \seq a=f(f(a))$ is not derivable in $EQ_M^-$. Since it is derivable in $EQ_M$, we have that  the contraction rule is not eliminable from derivations of $EQ_M$. As a consequence it is not eliminable in $EQ$ either.  In fact $a=f(a) \seq a=f(f(a))$ is derivable in $EQ$,  but it cannot have  a derivation in $EQ$ without applications of the contraction rule, since, by Proposition \ref{equiv},  such a derivation could be translated into a derivation in $EQ_M^-$ of  $a=f(a) \seq a=f(f(a))$ which we know it does not exist.

If we replace in $EQ_M$ or $EQ$ the cut rule by its context sharing version, namely the rule

\[
\ba{c}
\Gamma\seq A~~~\Gamma, A \seq H\\
\cline{1-1}
\Gamma \seq H
\ea
\]

then the contraction rule turns out to be derivable, thanks to the weakening and exchage rule, as shown by the following derivation:

\[
\ba{ccc}
F\seq F&&\\
\cline{1-1}
\Gamma, F\seq F&~~~&\Gamma, F,F \seq H\\
\cline{1-3}
\multicolumn{3}{c}{ \Gamma, F \seq H}
\ea
\]

In the case of $EQ$, a less trivial way of making the contraction rule eliminable is to replace the CNG rule by its context sharing version, 
while retaining the context independent cut rule (thus obtaining the system  denoted by $EQ$ in \cite{PP13})

\begin{definition}
$EQ'$ is obtained by replacing in $EQ$ the rule $CNG$ by its context sharing version $CNG'$, namely:

\[
\ba{c}
\Gamma \seq D\{v/r\}~~~~\Gamma \seq r=s\\
\cline{1-1}
\Gamma \seq D\{v/s\}
\ea
\]
$cf.EQ'$ is $EQ'$ deprived of the cut rule, and $ccf.EQ'$ is $cf.EQ'$ deprived also of the contraction rule.
\end{definition}



{\bf Notation} $\Gamma_F$ will denote  the  sequence  that is  obtained from $\Gamma$ by eliminating all the occurrences of $F$ but the last one, provided that there is at least one  occurrence of $F$  in $\Gamma$, and $\Gamma$ otherwise.

\begin{lemma}\label{cnf-ammissibilita'}
\bi
\item[a)] If $\Gamma \seq H$ is derivable in $ccf.EQ'$, then $\Gamma_F \seq H$ is  derivable in $ccf.EQ'$.

\item[b)] If $\Gamma, F, F \seq H$ is derivable in $ccf.EQ'$ then   $\Gamma, F \seq H$ is  derivable in $ccf.EQ'$.

\item[c)] If $\Gamma \seq H$ is derivable in $ccf.EQ'$ and $\Gamma_0$ contains all the formulae occurring in $\Gamma$ then 
$\Gamma_0\seq H$ is derivable in $ccf.EQ'$.

\ei

\end{lemma}

{\bf Proof}
$a)$ If $F$ has no occurrences in $\Gamma$ the conclusion is trivial. Otherwise we proceed
by induction on the height of a derivation ${\cal D}$ of  $\Gamma \seq H$ in $ccf.EQ$.
 In the base case
the given derivation reduces to the axiom $F\seq F$. Then    $\Gamma_F\seq H$ also reduces to the axiom $F \seq F$. As for the induction step, we have the following cases.
Case 1. ${\cal D}$ ends with an exchange. Then it suffices to apply the induction hypothesis to the immediate subderivation of ${\cal D}$, and then the same exchange by which ${\cal D}$ ends, unless one of the exchanded formula is $F$ itself, but it is not the last occurrence of $F$ in $\Gamma$. In that case the desired derivation is directly provided by the induction hypothesis.
Case 2. ${\cal D}$ 
ends with a weakening, i.e. is of the form:

\[
\ba{c}
{\cal D}_0\\
\Gamma'\seq H\\
\cline{1-1}
\Gamma', G\seq H
\ea
\]
By induction hypothesis there is a derivation of $\Gamma'_F \seq H$. If $F$ is different from $G$, then,
since $\Gamma_F$ coincides with $\Gamma'_F,G$, 
it suffices to apply  the same weakening to obtain the desired derivation of $\Gamma_F \seq H$. If $F$ coincides with $G$ and does not occur in $\Gamma'$
then ${\cal D}$ is already a derivation of $\Gamma_F\seq H$. Otherwise, if $F$ occurs last in $\Gamma'$ we are done. If not,  the desired derivation is obtained by applying the exchanges needed to bring the unique occurrence of $F$ in
$\Gamma'_F$ at the end of the sequence.

Case 3. ${\cal D}$ ends with a $CNG'$-inference. Since a $CNG'$-inference does not modify the antecedent of the premisses, the claim is an immediate consequence of the induction hypothesis.

$b)$ By $a)$,  if $\Gamma, F,F\seq H$ is derivable in $ccf.EQ$ and $\Gamma_0$ is obtained from $\Gamma$ by eliminating all the occurrences of $F$, then $\Gamma_0, F\seq H$ has a derivation  in $ccf.EQ$, from which by means of weakenings, introducing $F$, and exchanges we 
can obtain a derivation in $ccf.EQ$ of $\Gamma, F \seq H$.

$c)$ is obtained by applying  $a)$ for $k$ times, where $k$ is the number of different formulae occurring in $\Gamma$. $\Box$


\begin{proposition}
The contraction rule is  eliminable from derivations in $cf.EQ'$.
\end{proposition}

{\bf Proof}  By the previous Lemma \ref{cnf-ammissibilita'} $b)$, the contraction rule is admissible in $ccf.EQ'$
and therefore eliminable from derivations in $cf.EQ'$ $\Box$

\

Thus, taking into account the eliminability of the cut rule from derivations in  $EQ'$, established in \cite{PP13}, we have the following:

\begin{corollary}
Both the cut  and the contraction rule are   eliminable from derivations in $EQ'$.
\end{corollary}

To sum up: contraction elimination does not hold for $EQ_M$ and $EQ$, but it does hold if the cut rule, in the case of $EQ_M$, and the cut rule or the rule  $CNG$, in the case of $EQ$,  are replaced by their context sharing versions.

\end{document}